\documentclass[12pt]{article}
 
 \normalsize
 \usepackage{amssymb}		
 \usepackage{amstext}
 \usepackage{amsfonts}
 \usepackage{amsmath}		
 \usepackage{latexsym}		
 \usepackage{epsfig,epic}       

 \newtheorem{definition}{Definition}[section]
 \newtheorem{conjecture}[definition]{Conjecture}
 \newtheorem{corollary}[definition]{Corollary}
 \newtheorem{example}[definition]{Example}
 \newtheorem{lemma}[definition]{Lemma}
 \newtheorem{proposition}[definition]{Proposition}
 \newtheorem{remark}[definition]{Remark}
 
 \newtheorem{theorem}[definition]{Theorem}
 \newtheorem{warning}[definition]{Warning}

 \newcommand{\para}{\medskip \indent}
 \newcommand{\proof}{\noindent\textsc{Proof.\ }}
 \newcommand{\discussproof}{\noindent\textsc{Discussion of proof.\ }}
 \newcommand{\qed}{\hspace*{\fill}\(\Box\) \medskip}
 \newcommand{\one}{\ensuremath{(\mathrm{i})}}
 \newcommand{\two}{\ensuremath{(\mathrm{ii})}}
 \newcommand{\three}{\ensuremath{(\mathrm{iii})}}

 \newcommand{\C}{\ensuremath{\mathbb{C}}}
 
 \newcommand{\Q}{\ensuremath{\mathbb{Q}}}
 \newcommand{\R}{\ensuremath{\mathbb{R}}}
 \newcommand{\Z}{\ensuremath{\mathbb{Z}}}
 
 \newcommand{\U}{\ensuremath{\operatorname{U}}}
 \newcommand{\SU}{\ensuremath{\operatorname{SU}}}
 \newcommand{\SL}{\ensuremath{\operatorname{SL}}}
 \newcommand{\GL}{\ensuremath{\operatorname{GL}}}
 \newcommand{\oo}{\ensuremath{\infty}}

 \newcommand{\Conj}{\ensuremath{\operatorname{Conj}}} 
 \newcommand{\Gr}{\ensuremath{\operatorname{Gr}}}
 \newcommand{\Hom}{\ensuremath{\operatorname{Hom}}}

 \newcommand{\Spec}{\ensuremath{\operatorname{Spec}}}
 
 \newcommand{\Var}{\ensuremath{\mathcal{V}_{\mathbb{C}}}}

 \newcommand{\age}{\operatorname{age}}
 
 \newcommand{\blank}{\operatorname{\;\;\;}}

 \newcommand{\dbrk}[1]{\left[\kern-.15em\left[#1\right]\kern-.15em\right]}
 
 \newcommand{\diag}{\ensuremath{\operatorname{diag}}}
 \newcommand{\dimn}{\operatorname{dim}}
 \newcommand{\dmu}{\ensuremath{\operatorname{d\mu}}}

 \newcommand{\lm}{\ensuremath{\mathbb{L}}}

 \title{An introduction to motivic integration}
 \author{Alastair Craw}

\begin{document}
 
 \maketitle
 
 \begin{abstract}
 By associating a `motivic integral' to every complex projective variety \(X\) with at worst canonical,  Gorenstein singularities,  Kontsevich~\cite{Kontsevich:lao} proved that,  when there exists a crepant resolution of singularities \(\varphi\colon Y \rightarrow X\),  the Hodge numbers of \(Y\) do not depend upon the choice of the resolution.  

 In this article we provide an elementary introduction to the theory of motivic integration,  leading to a proof of the result described above.  We calculate the motivic integral of several quotient singularities and discuss these calculations in the context of the cohomological McKay correspondence.
 \end{abstract}
 
 \section{Introduction}
 \label{sec:intro}
 String theory is the motivation behind the `Mirror Symmetry Conjecture' which states that there exist mirror pairs of Calabi--Yau varieties with `certain compatibilities'.  For instance,  if \( (X,X^*)\) is a smooth, projective mirror pair of dimension \(n\),  then we expect the relation
 \begin{equation}
 \label{eqn:smoothmirrorpair}
h^{p,q}(X) = h^{n-p,q}(X^*)
 \end{equation}
 to hold between their Hodge numbers.  Mirror pairs are not smooth in general and the compatibility relation (\ref{eqn:smoothmirrorpair}) can fail to hold if either \(X\) or its mirror have singularities.  In this case,  string theory led to the following revised compatibility relation:  if there exist crepant resolutions $Y \rightarrow X$ and \(Y^* \rightarrow X^*\) then we expect the relation
 \begin{equation}
 \label{eqn:smoothmirrorpair2}
h^{p,q}(Y) = h^{n-p,q}(Y^*)
 \end{equation}
 to hold between the Hodge numbers of the smooth varieties \(Y\) and \(Y^*\).  A resolution \(\varphi\colon Y \rightarrow X\) is said to be crepant if \(K_Y = \varphi^*X\).

 \indent However,  it is not obvious that the revised relation (\ref{eqn:smoothmirrorpair2}) is well defined; if a crepant resolution exists it is not necessarily unique.  In particular,  given two crepant resolutions \(Y_1 \rightarrow X\) and \(Y_2 \rightarrow X\), it is not clear a priori that the Hodge numbers of \(Y_1\) and \(Y_2\) are equal.

 \indent Nevertheless,  the consistency of string theory led Batyrev and Dais~\cite{Batyrev:smc} to conjecture that,  when \(X\) has only mild Gorenstein singularities,  the Hodge numbers of \(Y_1\) and \(Y_2\) are equal.  In a subsequent paper~\cite{Batyrev:bcy},  Batyrev used methods of \(p\)-adic integration to prove that the Betti numbers of \(Y_1\) and \(Y_2\) are equal.  Kontsevich~\cite{Kontsevich:lao} later proved that the Hodge numbers are equal by introducing the notion of motivic integration. 

 \bigskip

 This article provides an elementary introduction to Kontsevich's theory of motivic integration.  The first step is to construct the motivic integral of a pair \( (Y,D)\),  for a complex manifold \(Y\) and an effective divisor \(D\) on \(Y\) with simple normal crossings.  We define the \emph{space of formal arcs} \(J_\infty(Y)\) of \(Y\) and associate a function \(F_D\) defined on \(J_\infty(Y)\) to the divisor \(D\).  The \emph{motivic integral} of the pair \( (Y,D)\) is the integral of \(F_D\) over \(J_\infty(Y)\) with respect to a certain measure \(\mu\) on \(J_\infty(Y)\).  This measure is not real-valued;  the subtlety in the construction is in defining the ring in which \(\mu\) takes values.  We adopt the structure of the proof of Theorem~6.28 from Batyrev~\cite{Batyrev:shn} to establish the following user-friendly formula:

 \begin{theorem}[formula for the motivic integral]
 \label{introthm:motintformula}
 Let \(Y\) be a complex manifold of dimension \(n\) and \(D = \sum_{i=1}^r a_iD_i\) an effective divisor on \(Y\) with simple normal crossings.  The motivic integral of the pair \( (Y,D)\) is 
 \begin{equation}
 \label{introeqn:motintformula}
 \int_{J_\infty(Y)} F_{D} \dmu  =  \sum_{J\subseteq\{1,\dots ,r\} }[D_J^\circ]\cdot \left(\prod_{j\in J}\frac{\lm - 1}{\lm^{a_j + 1} - 1}\right)\cdot \lm^{-n}
 \end{equation}
 where we sum over all subsets \(J\subseteq\{1,\dots ,r\}\) including \(J = \emptyset\).
 \end{theorem}

 \indent The motivic integral of a complex algebraic variety \(X\) with Gorenstein canonical singularities is defined to be the motivic integral of a pair \( (Y,D)\),  where \(Y\rightarrow X\) is a resolution of singularities for which the discrepancy divisor \(D\) has simple normal crossings.  Crucially,  this is well defined independent of the choice of resolution.  The motivic integral induces a \emph{stringy \(E\)-function} 
 \begin{equation}
 \label{introeqn:stringyE}
 E_{\mathrm{st}}(X) := \sum_{J\subseteq\{1,\dots ,r\} }E(D_J^\circ)\cdot\left(\prod_{j\in J}\frac{uv - 1}{(uv)^{a_j + 1} - 1}\right)
 \end{equation}
 which is also independent of the choice of resolution (see Warning~\ref{warning:error}).  The \(E\)-polynomials \(E(D_J^\circ)\) encode the Hodge--Deligne numbers of open strata \(D_J^\circ\subset Y\),  and the stringy \(E\)-function records these numbers with certain `correction terms' written in parentheses in formula (\ref{introeqn:stringyE}).  When \(Y\to X\) is a crepant resolution the correction terms disappear leaving simply the terms \(E(D_J^\circ)\) whose sum is the \(E\)-polynomial of \(Y\).  In this way the function \(E_{\mathrm{st}}(X)\) encodes the Hodge numbers of a crepant resolution \(Y\to X\),  thereby establishing Kontsevich's result on the equality of Hodge numbers (see Section~\ref{sec:mni} for a brief look at the motivic nature of the integral).

 However,  crepant resolutions do not exist in general.  To get a better feeling for the stringy \(E\)-function of varieties admiting no crepant resolution we calculate \(E_{\mathrm{st}}(X)\) for several 4- and 6-dimensional Gorenstein terminal cyclic quotient singularities.  These calculations lead naturally to a discussion of how Batyrev used the theory of motivic integration to prove a refined version of the generalised McKay correspondence conjecture of Reid~\cite{Reid:pen}:

 \begin{theorem}[Batyrev]
 \label{introthm:gmc}
 Let \(G\subset \SL(n,\C)\) be a finite subgroup and suppose that the quotient \(X = \C^{n}/G\) admits a crepant resolution \(\varphi\colon Y\rightarrow X\).  Then \(H^{*}(Y,\Q)\) has a basis consisting of algebraic cycles corresponding one-to-one with conjugacy classes of \(G\).  In particular,  the Euler number of \(Y\) equals the number of conjugacy classes of \(G\).
 \end{theorem}

 This result generalises a theorem of McKay~\cite{McKay:gsfg}, namely,  that the graph of ADE type associated to a Kleinian singularity \(\C^{2}/G\) can be constructed using only the representation theory of the finite subgroup \(G \subset \SL(2,\C)\).  This establishes a one-to-one correspondence between a basis for the cohomology of the minimal resolution \(Y\) of \(\C^{2}/G\) and the irreducible representations of \(G\),  hence equality between the Euler number \(e(Y)\) and the number of irreducible representations (or conjugacy classes) of \(G\). 

 Batyrev~\cite{Batyrev:nai,Batyrev:ca} derives Theorem~\ref{introthm:gmc} as a consequence of the equality between the stringy \(E\)-function of the quotient \(\C^{n}/G\) and the `orbifold \(E\)-function' of the pair \((\C^{n},G)\).   We choose not to discuss the orbifold function here;  instead, we present a simple,  direct proof of Theorem~\ref{introthm:gmc} for a finite Abelian subgroup \(G\subset \SL(n,\C)\).





 \indent The original references on this topic are Batyrev~\cite[\S 6]{Batyrev:shn} and Denef and Loeser~\cite{Denef:goa}.  The more recent article by Looijenga~\cite{Looijenga:mm} provides a detailed survey of motivic integration.

\bigskip

\noindent \textbf{Acknowledgements\ } I wish to thank Miles Reid for his encouragement,  comments and corrections,  and Willem Veys who provided me with a counterexample to a conjecture which appeared in an earlier version of this article.  Thanks also to Victor Batyrev and Alastair King for their comments during my recent PhD thesis defence.

 \section{Construction of the motivic integral}
 \label{sec:cmi}

 \subsection{The space of formal arcs of a complex manifold}
 \label{sec:sasv}
 \begin{definition} 
 \emph{Let \(Y\) be a complex manifold of dimension \(n\),  and \(y\in Y\) a point.  A \emph{\(k\)-jet} over \(y\) is a morphism
 \[
 \gamma_y\colon\Spec\ \mathbb{C}[z]/\langle z^{k+1}\rangle \longrightarrow Y\quad \mbox{with}\quad \gamma_y(\Spec\ \C) = y.
 \]
Once local co-ordinates are chosen the space of \(k\)-jets over \(y\) can be viewed as the space of \(n\)-tuples of polynomials of degree \(k\) whose constant terms are zero.  Let \(J_{k}(Y)\) denote the bundle over \(Y\) whose fibre over \(y\in Y\) is the space of \(k\)-jets over \(y\).  A \emph{formal arc} over \(y\) is a morphism
 \[
 \gamma_y\colon\Spec\ \mathbb{C}\dbrk{z} \longrightarrow Y\quad \mbox{with}\quad \gamma_y(\Spec\ \C) = y.
 \]
 Once local co-ordinates are chosen the space of formal arcs over \(y\) can be viewed as the space of \(n\)-tuples of power series whose constant terms are zero.  Let \(\pi_{0}\colon J_{\oo}(Y) \rightarrow Y\) denote the bundle whose fibre over \(y\in Y\) is the space of formal arcs over \(y\).  For each \(k\in \mathbb{Z}_{\geq 0}\) the inclusion \(\mathbb{C}[z]/\langle z^{k+1}\rangle \hookrightarrow \C\dbrk{z}\) induces a surjective map 
 \[
 \pi_{k}\colon J_{\oo}(Y) \longrightarrow J_{k}(Y).
 \]
 } 
 \end{definition}

\begin{definition}\label{defn:cylinder} \emph{A subset \(C\subseteq J_\infty(Y)\) of the space of formal arcs is called a \emph{cylinder set} if \(C = \pi_k^{-1}(B_k)\) for \(k\in \mathbb{Z}_{\geq 0}\) and \(B_k \subseteq J_k(Y)\) a constructible subset.  Recall that a subset of a variety is \emph{constructible} if it is a finite, disjoint union of (Zariski) locally closed subvarieties.}\end{definition}

\para It's clear that the collection of cylinder sets forms an algebra of sets (see \cite[p.\ 10]{Rudin:rca});  that is,  \(J_\oo(Y) = \pi_0^{-1}(Y)\) is a cylinder set,  as are finite unions and complements (and hence finite intersections) of cylinder sets.

\subsection{The function \protect\(F_D\protect\) associated to an effective divisor}
\label{sec:fad}
\begin{definition} \emph{Let \(D\) be an effective divisor on \(Y\),  \(y\in Y\) a point,  and \(g\) a local defining equation for \(D\) on a neighbourhood \(U\) of \(y\).  For an arc \(\gamma_u\) over a point \(u\in U\),  define the \emph{intersection number} \(\gamma_u\cdot D\) to be the order of vanishing of the formal power series \(g(\gamma_u(z))\) at \(z = 0\).  The function 
\[F_D \colon J_\infty(Y) \rightarrow \mathbb{Z}_{\geq 0}\cup \oo\] 
associated to the divisor \(D\) on \(Y\) is given by \(F_D(\gamma_u) = \gamma_u\cdot D\).  If we write \(D = \sum_{i=1}^r a_iD_i\) as a linear combination of prime divisors then \(g\) decomposes as a product \(g = \textstyle{\prod_{i=1}^r} g_i^{a_i}\) of defining equations for \(D_i\),  hence \(F_{D} = \sum_{i=1}^r a_{i}F_{D_i}\).  Furthermore
\begin{equation}
F_{D_i}(\gamma_u) = 0 \iff u\notin D_i \quad \mbox{and}\quad F_{D_i}(\gamma_u) = \oo \iff \gamma_u \subseteq D_i. \label{eqn:zero}
\end{equation}
}\end{definition}

\para Our ultimate goal is to integrate the function \(F_D\) over \( J_\oo(Y)\),  so we must understand the nature of the level set \(F_{D}^{-1}(s) \subseteq J_\oo(Y)\) for each \(s \in \mathbb{Z}_{\geq 0} \cup \oo\).  With this goal in mind,  we introduce a partition of \(F_D^{-1}(s)\).

\begin{definition}
\label{defn:partition_level_set}
\emph{For \(D = \sum_{i=1}^r a_iD_i\) and \(J\subseteq \{1,\dots ,r\}\) any subset,  define
\[D_J := \left\{\begin{array}{cc}\bigcap_{j\in J}D_j & \mbox{if\ }J\neq \emptyset \\ Y &  \mbox{if\ }J = \emptyset \end{array}\right. \hspace{5mm}\mbox{and}\hspace{5mm} D_J^\circ := D_J \setminus \bigcup_{i\in \{1,\dots ,r\}\setminus J} D_i.\]
These subvarieties stratify \(Y\) and define a partition of the space of arcs into cylinder sets:
\[Y = \bigsqcup_{J\subseteq \{1,\dots ,r\}} D^\circ _J\quad\mbox{and}\quad J_\oo(Y) = \bigsqcup_{J\subseteq \{1,\dots ,r\}} \pi_0^{-1}(D^\circ _J).\]
For any \(s\in \mathbb{Z}_{\geq 0}\) and \(J\subseteq \{1,\dots ,r\}\),  define
\[M_{J,s} :=  \left\{(m_1,\dots ,m_r) \in \mathbb{Z}^r_{\geq 0}\bigm{|}  \sum a_im_i = s\;\mbox{with\ } m_j > 0 \Leftrightarrow j\in J\right\}.\]
It now follows from (\ref{eqn:zero}) that
\[\gamma_u \in \pi_0^{-1}(D_J^\circ)\cap F_D^{-1}(s)  \iff \left(F_{D_1}(\gamma_u),\dots ,F_{D_r}(\gamma_u)\right) \in M_{J,s}.\]
As a result we produce a finite partition of the level set
\begin{equation}
\label{eqn:F_D_partition}
F_{D}^{-1}(s) = \bigsqcup_{J\subseteq \{1,\dots ,r\}}\bigsqcup_{(m_1,\dots ,m_r) \in M_{J,s}} \left(\textstyle{\bigcap_{i = 1,\dots r}}F_{D_i}^{-1}(m_i)\right).
\end{equation}
}\end{definition}

\begin{proposition} \label{prop:cyl} If \(D\) is an effective divisor with simple normal crossings then \(F_D^{-1}(s)\) is a cylinder set (see Definition \ref{defn:cylinder}) for each \(s \in \mathbb{Z}_{\geq 0}\).\end{proposition}

Recall (see \cite[p.\ 25]{Kawamata:mmp}) that a divisor \(D = \sum_{i=1}^r a_iD_i\) on \(Y\) has only simple normal crossings if,  at each point \(y\in Y\),  there is a neighbourhood \(U\) of \(y\) with coordinates \(z_1,\dots ,z_n\) for which a local defining equation for \(D\) is 
\begin{equation}
\label{eqn:snc_equation}
g = z_1^{a_1}\cdots z_{j_y}^{a_{j_y}}\quad\mbox{for some\ } j_y \leq n.
\end{equation}  

\medskip

\noindent\textsc{Proof of Proposition~\ref{prop:cyl}.\ } A finite union of cylinder sets is cylinder and we have a partition (\ref{eqn:F_D_partition}) of  \(F_D^{-1}(s)\),  so it is enough to prove,  for some \(J\subseteq \{1,\dots ,r\}\) and \((m_1,\dots ,m_r) \in M_{J,s}\),  that \(\textstyle{\bigcap_{i = 1,\dots r}}F_{D_i}^{-1}(m_i)\) is a cylinder set\footnote{Finite intersections of cylinder sets are cylinder,  so we could reduce to proving the result for \(F_{D_i}^{-1}(m_i)\).  However we require (\ref{eqn:U_constructible}) in \S\ref{sec:mioYD}.}.  Cover \(Y = \bigcup U\) by finitely many charts on which \(D\) has a local equation of the form (\ref{eqn:snc_equation}),  and lift to cover \(J_\infty (Y) = \bigcup \pi_0^{-1}(U)\).  Clearly we need only prove that the set
\[U_{m_1,\dots,m_r}:= \textstyle{\bigcap_{i = 1,\dots r}}F_{D_i}^{-1}(m_i)\cap \pi_0^{-1}(U)\] 
is cylinder.  In the notation of (\ref{eqn:snc_equation}),  if \(J \nsubseteq \{1,\dots ,j_y\}\) then \(D_J^\circ \cap U = \emptyset\) which forces \(U_{m_1,\dots,m_r}\subset \pi_0^{-1}(D_J^\circ \cap U)\) to be empty,  and hence a cylinder set.  We suppose therefore that \(J \subseteq \{1,\dots ,j_y\}\),  thus \(\vert J\vert \leq n\) holds by (\ref{eqn:snc_equation}).

\para The key observation is that when we regard each arc \(\gamma_u\) as an \(n\)-tuple \( (p_1(z),\dots ,p_n(z) )\) of formal power series with zero constant term,  each condition \(F_{D_i}(\gamma_u) = m_i\) is equivalent to a condition on the truncation of the power series \(p_i(z)\) to degree \(m_i\).  Indeed,  since \(D_i\) is cut out by \( z_i = 0\) on \(U\),  it follows that \(F_{D_i}(\gamma_u) = \{\mbox{order of\ }p_i(z)\; \mbox{at\ } z = 0\}\).  Thus \(\gamma_u \in F_{D_i}^{-1}(m_i)\) if and only if the truncation of \(p_i(z)\) to degree \(m_i\) is of the form \(c_{m_i} z^{m_i}\),  with \(c_{m_i}\neq 0\).  Truncating all \(n\) of the power series to degree \(t := \max\{m_j | j\in J\}\) produces \(n - \vert J\vert\) polynomials of degree \(t\) with zero constant term,  and,  for each \(j\in J\),  a polynomial of the form 
\[\pi_t(p_j(z)) = 0 + \cdots + 0 + c_{m_j}z^{m_j} + c_{(m_j+1)}z^{m_j + 1} + \dots + c_{t}z^t\]
for \(c_{m_j}\in \C^*\) and \(c_k \in \C \;\: \forall \:k > m_j\).  The space of all such \(n\)-tuples is isomorphic to \(\C^{t(n-\vert J\vert)}\times(\C^*)^{\vert J \vert}\times\C^{t\vert J\vert - \sum_{j\in J}m_j}\), hence
\begin{equation}
\label{eqn:U_constructible}
U_{m_1,\dots,m_r} = \pi_t^{-1}\left((U\cap D_J^\circ)\times \C^{tn - \sum_{j\in J}m_j}\times(\C^*)^{\vert J \vert}\right).
\end{equation}
The set \( (U\cap D_J^\circ)\times \C^{tn - \sum_{j\in J}m_j }\times(\C^*)^{\vert J \vert}\) is constructible,  so \(U_{m_1,\dots,m_r}\) is a cylinder set.  This completes the proof of the proposition.\qed

\para It is worth noting that \(F_D^{-1}(\oo)\) is not a cylinder set. Indeed,  suppose otherwise,  so there exists a constructible subset \(B_k \subseteq J_k(Y)\) for which \(F_D^{-1}(\oo) = \pi_k^{-1}(B_k)\).  Each arc \(\gamma_y \in F_D^{-1}(\oo)\) is an \(n\)-tuple of power series,  at least one of which is identically zero,  whereas each \(\gamma_y \in \pi_k^{-1}(B_k)\) is an \(n\)-tuple of power series whose terms of degree higher than \(k\) may take any complex value;  clearly this is absurd.

\begin{proposition} \label{prop:notcyl} \(F_D^{-1}(\oo)\) is a countable intersection of cylinder sets.\end{proposition}
\proof Observe that 
\begin{equation}
\label{eqn:countint}
F_D^{-1}(\oo) = \bigcap_{k\in \mathbb{Z}_{\geq 0}} \pi_k^{-1} \pi_k\big{(}F_{D}^{-1}(\infty)\big{)}
\end{equation}
because a power series is identically zero if and only if its truncation to degree \(k\) is the zero polynomial,  for all \(k\in \mathbb{Z}_{\geq 0}\).  It is easy to see that the sets \(\pi_k(F_{D}^{-1}(\infty))\subset J_k(Y)\) are constructible.  \qed

\subsection{A measure \protect\(\mu\protect\) on the space of formal arcs}
\label{sec:msa}
In this section we define a measure \(\mu\) on \(J_\oo(Y)\) with respect to which the function \(F_D\) is measurable.  The measure is not real-valued,  so we begin this section by constructing the ring in which \(\mu\) takes values.

 \begin{definition}
 \label{defn:grav} 
 \emph{Let \Var\ denote the category of complex algebraic varieties.  The Grothendieck group of \Var\ is the free Abelian group on the isomorphism classes \( [V]\) of complex algebraic varieties modulo the subgroup generated by elements of the form \( [V] - [V'] - [V\: \setminus V']\) for a closed subset \(V'\subseteq V\).  The product of varieties induces a ring structure \( [V] \cdot [V'] = [V\times V']\),  and the resulting ring,  denoted by \(K_0(\Var)\),  is called the \emph{Grothendieck ring of complex algebraic varieties}.  Let} 
 \[
 [ \blank ] \colon \mbox{\emph{Ob}}\Var \longrightarrow K_0(\Var)
 \]
 \emph{denote the natural map sending \(V\) to its class \([V]\) in the Grothendieck ring.  This map is universal with respect to maps which are additive on disjoint unions of constructible subsets,  and which respect products.
 }
 \end{definition} 
 
 \indent  Write\footnote{See the Appendix~\ref{sec:mni}:  the class of \C\ in \(K_0(\Var)\) corresponds to the Tate motive \(\lm\).} \(1:= [\mbox{point}]\) and \(\lm:= [\C]\).  Then
\[ [\C^*] = [\C -\{0\}] = [\C] - [\{0\}] = \lm - 1.\]
Also,  if \(f\colon Y \rightarrow X\) is a locally trivial fibration w.r.t.\ the Zariski topology and \(F\) is the fibre over a closed point then \([Y] = [F\times X]\).

\begin{definition} \emph{Let \(K_0(\Var)[\lm^{-1}] := S^{-1} K_0(\Var)\) denote the ring of fractions of \(K_0(\Var)\) with respect to the multiplicative set \(S := \{1,\lm,\lm^2,\dots \}\).}\end{definition}


\begin{definition} \label{defn:mutilde} \emph{Recall that cylinder sets in \(J_\oo(Y)\) are subsets \(\pi_k^{-1}(B_k)\subset J_\infty(Y)\) for \(k\in \mathbb{Z}_{\geq 0}\) and for \(B_k \subseteq J_k(Y)\) a constructible subset.  The function 
\[\widetilde{\mu}\colon  \Big{\{}\mbox{cylinder sets in\ }J_\oo(Y)\Big{\}}\; \longrightarrow  K_0(\Var)[\lm^{-1}] \]
 which assigns a `measure' to each cylinder set is defined by
\[\pi_k^{-1}(B_k)\rightarrow [B_k]\cdot \lm^{-n(k+1)}.\]
Using the fact that the map \( [ \;\;\; ]\) introduced in Definition~\ref{defn:grav} is additive on disjoint unions of constructible sets,  it is straightforward to show that
\[\widetilde{\mu} \left(\bigsqcup_{i=1}^l C_i\right) = \sum_{i=1}^l \widetilde{\mu}(C_i)\quad\text{for cylinder sets\ } C_1,\dots ,C_l.\]
For this reason we call \(\widetilde{\mu}\) a finitely additive measure.}\end{definition}

\para Proposition~\ref{prop:cyl} states that for \(s\in \mathbb{Z}_{\geq 0}\),  the level set \(F_D^{-1}(s)\) is a cylinder set,  and is therefore \(\widetilde{\mu}\)-measurable.  However,  \(F_D\) is not \(\widetilde{\mu}\)-measurable because \(F_D^{-1}(\oo)\) is not cylinder.  To proceed,  we extend \(\widetilde{\mu}\) to a measure \(\mu\) with respect to which \(F_D^{-1}(\oo)\) is measurable.  

 \indent The following discussion is intended to motivate the definition of \(\mu\) (see Definition~\ref{defn:mu} to follow).  The set \(J_\oo(Y) \setminus F_D^{-1}(\oo)\) is a countable disjoint union of cylinder sets
\begin{equation}
J_\oo(Y) \setminus \pi_0^{-1} \pi_0(F_{D}^{-1}(\infty)) \sqcup \bigsqcup_{k\in \mathbb{Z}_{\geq 0}} \bigl{(}\pi_k^{-1} \pi_k(F_{D}^{-1}(\infty))\setminus\pi_{k+1}^{-1} \pi_{k+1}(F_{D}^{-1}(\infty))\bigr{)}; \label{eqn:countunion}
\end{equation}
to see this,  take complements in equation (\ref{eqn:countint}) of Proposition \ref{prop:notcyl}.  Our goal is to extend \(\widetilde{\mu}\) to a measure \(\mu\) defined on the collection of countable disjoint unions of cylinder sets so that the set \(J_\oo(Y) \setminus F_D^{-1}(\oo)\),  and hence its complement \(F_D^{-1}(\oo)\), is \(\mu\)-measurable.  One would like to define
\begin{equation}
\label{eqn:mu}\mu\left(\bigsqcup_{i\in\mathbb{N}} C_i\right) := \sum_{i\in\mathbb{N}}\mu(C_i) = \sum_{i\in\mathbb{N}} \widetilde{\mu}(C_i)\quad\text{for cylinder sets\ } C_1,\dots ,C_l.\end{equation}

\para However,  countable sums are not defined in \(K_0(\Var)[\lm^{-1}]\).  Furthermore,  given a countable disjoint union \(C = \bigsqcup_{i\in\mathbb{N}} C_i\),  it is not clear a priori that \(\mu(C)\) defined by formula (\ref{eqn:mu}) is independent of the choice of the \(C_i\).  

\indent Kontsevich~\cite{Kontsevich:lao} solved both of these problems at once by completing the ring \(K_0(\Var)[\lm^{-1}]\),  thereby allowing appropriate countable sums,  in such a way that the measure of the set \(C = \bigsqcup_{i\in\mathbb{N}} C_i\) is independent of the choice of the \(C_i\),  assuming that \(\mu(C_i)\rightarrow 0\) as \(i \rightarrow \oo\). 

\begin{definition}\label{defn:R} \emph{Let \(R\) denote the completion of the ring \(K_0(\Var)[\lm^{-1}]\) with respect to the filtration 
\[\dots \supseteq F^{-1}K_0(\Var)[\lm^{-1}] \supseteq F^0K_0(\Var)[\lm^{-1}]\supseteq F^1K_0(\Var)[\lm^{-1}]\supseteq \cdots \]
where for each \(m\in \mathbb{Z}\),  \(F^mK_0(\Var)[\lm^{-1}]\) is the subgroup of \(K_0(\Var)[\lm^{-1}]\) generated by elements of the form \([V]\cdot\lm^{-i}\) for \(i - \mbox{dim\ }V \geq m\).  The natural completion map is denoted \(\phi\colon K_0(\Var)[\lm^{-1}] \longrightarrow R\).}\end{definition}

\para By composing \(\widetilde{\mu}\) with the natural completion map \(\phi\),  we produce a finitely additive measure with values in the ring \(R\),  namely 
\[\widetilde{\mu} := \phi \circ\widetilde{\mu}\colon \pi_k^{-1}(B_k) \rightarrow \phi\left([B_k]\cdot \lm^{-n(k+1)}\right)\]
which we also denote \(\widetilde{\mu}\).  Given a sequence of cylinder sets \(\{C_i\}\) one may now ask whether or not \(\widetilde{\mu}(C_i)\rightarrow 0\) as \(i \rightarrow \oo\).  We are finally in a position to define the measure \(\mu\) on the space of formal arcs.

\begin{definition} \label{defn:mu}\emph{Let \(\mathcal{C}\) denote the collection of countable disjoint unions of cylinder sets \(\textstyle{\bigsqcup_{i\in \mathbb{N}}} C_i\) for which \(\widetilde{\mu}(C_i) \rightarrow  0\) as \(i \rightarrow \oo\),  together with the complements of such sets.  Extend \(\widetilde{\mu}\) to a measure \(\mu\) defined on \(\mathcal{C}\) which takes values in \(R\) given by
\[\bigsqcup_{i\in \mathbb{N}} C_i\longrightarrow \sum_{i\in \mathbb{N}} \widetilde{\mu}(C_i).\]
It is nontrivial to show (see \cite[\S 3.2]{Denef:goa} or \cite[\S 6.18]{Batyrev:shn}) that this definition is independent of the choice of the \(C_i\).}
\end{definition}

\begin{proposition} \label{prop:FDmble} \(F_D\) is \(\mu\)-measurable,  and \(\mu(F_D^{-1}(\oo)) = 0\).\end{proposition}
\proof We prove that \(F_D^{-1}(\oo)\) (in fact its complement) lies in \(\mathcal{C}\).  It's clear from (\ref{eqn:countunion}) that we need only prove that \(\mu(\pi_k^{-1}\pi_k(F_{D}^{-1}(\infty)))\rightarrow 0\) as \(k\rightarrow \oo\).  Lemma~\ref{lemma:m'ble} below reveals that \(\mu(\pi_k^{-1}\pi_k(F_{D}^{-1}(\infty))) \in \phi(F^{k+1}K_0(\Var)[\lm^{-1}])\) which, by the nature of the topology on \(R\),  tends to zero as \(k\) tends to infinity.  This proves the first statement.  Using (\ref{eqn:countunion}) we calculate 
\begin{multline}
\mu\bigl{(}J_\oo(Y) \setminus F_D^{-1}(\oo)\bigr{)} = \widetilde{\mu}\bigl{(}J_\oo(Y) \setminus \pi_0^{-1} \pi_0(F_{D}^{-1}(\infty))\bigr{)} \\
 + \sum_{k\in \mathbb{Z}_{\geq 0}} \widetilde{\mu}\left(\pi_k^{-1} \pi_k(F_{D}^{-1}(\infty))\setminus\pi_{k+1}^{-1} \pi_{k+1}(F_{D}^{-1}(\infty))\right).
\end{multline}
This equals \(\mu(J_\oo(Y)) - \lim_{k\rightarrow \oo} \mu\left(\pi_k^{-1} \pi_k(F_{D}^{-1}(\infty))\right)\).  By the above remark,  this is simply \(\mu(J_\oo(Y))\),  so \(\mu(F_D^{-1}(\oo)) = 0\) as required.\qed

\begin{lemma} 
\label{lemma:m'ble}
\(\widetilde{\mu}(\pi_k^{-1}\pi_k(F_{D}^{-1}(\infty)))\in F^{k+1}K_0(\Var)[\lm^{-1}]\)
\end{lemma}
\proof It is enough to prove the result for a prime divisor \(D\),  since \(F_D^{-1}(\oo)\) is the union of sets \(F_{D_i}^{-1}(\oo)\).  Choose coordinates on a chart \(U\) in which \(D\) is \( (z_1 = 0)\).  Each \(\gamma_y \in F_{D}^{-1}(\infty)\cap \pi_0^{-1}(U)\) is an \(n\)-tuple \( (p_1(z),\dots ,p_n(z) )\) of power series over \(y\in U\cap D\) such that \(p_1(z)\) is identically zero.  Truncating these power series to degree \(k\) leaves \(n-1\) polynomials of degree \(k\) with zero constant term,  and the zero polynomial \(\pi_k(p_1(z))\).  The space of all such polynomials is isomorphic to \(\C^{(n-1)k}\),  so that \(\pi_k(F_{D}^{-1}(\infty)\cap \pi_0^{-1}(U))\simeq (U\cap D) \times \C^{(n-1)k}\).  Thus \([\pi_k(F_{D}^{-1}(\infty)] = [D]\cdot[\C^{(n-1)k}]\) and
\begin{eqnarray*}
\widetilde{\mu}(\pi_k^{-1}\pi_k(F_{D}^{-1}(\infty)) & = & [\pi_k(F_{D}^{-1}(\infty)]\cdot \lm^{-n(k+1)} \\
  & = & [D]\cdot \lm^{(n-1)k}\cdot \lm^{-n(k+1)} \\
  & = & [D]\cdot \lm^{-(n+k)} 
\end{eqnarray*}
This lies in \(F^{k+1}K_0(\Var)[\lm^{-1}]\) since \(D\) has dimension \(n-1\).   \qed

\subsection{The motivic integral of a pair \protect\( (Y,D)\protect\)}
\label{sec:mioYD}
\begin{definition}\label{defn:mi} \emph{Let \(Y\) be a complex manifold of dimension \(n\),  and choose an effective divisor \(D = \sum_{i=1}^r a_iD_i\) on \(Y\) with only simple normal crossings.  The \emph{motivic integral of the pair} \( (Y,D)\) is 
\[\int_{J_\infty(Y)} F_D \dmu := \sum_{s\in \mathbb{Z}_{\geq 0}\cup \infty} \mu\left(F_D^{-1}(s)\right) \cdot \lm^{-s}.\]
Since the set \(F_D^{-1}(\infty)\subset J_\oo(Y)\) has measure zero (see Proposition~\ref{prop:FDmble}),  we need only integrate over \(J_{\oo}(Y) \setminus F_D^{-1}(\infty)\),  so we need only sum over \(s\in \mathbb{Z}_{\geq 0}\).}\end{definition}

 \indent We now show that the motivic integral converges in the ring \(R\) introduced in Definition~\ref{defn:R}.  In doing so,  we establish a user-friendly formula.
    
\begin{theorem} \label{thm:formula} Let \(Y\) be a complex manifold of dimension \(n\) and \(D = \sum_{i=1}^r a_iD_i\) an effective divisor on \(Y\) with only simple normal crossings.  The motivic integral of the pair \( (Y,D)\) is 
\[\int_{J_\infty(Y)} F_D \dmu  =  \sum_{J\subseteq\{1,\dots ,r\} }[D_J^\circ]\cdot \left(\prod_{j\in J}\frac{\lm - 1}{\lm^{a_j + 1} - 1}\right)\cdot \lm^{-n}\]
where we sum over all subsets \(J\subseteq\{1,\dots ,r\}\) including \(J = \emptyset\).
\end{theorem}
\proof
In the proof of Proposition~\ref{prop:cyl} we cover \(Y\) by sets \(\{U\}\) and prove that \(\textstyle{\bigcap_{i = 1,\dots r}}F_{D_i}^{-1}(m_i)\cap \pi_0^{-1}(U)\) is a cylinder set of the form 
\[\pi_t^{-1}\left((U\cap D_J^\circ)\times \C^{tn - \sum_{j\in J}m_j}\times(\C^*)^{\vert J \vert}\right).\]
Since the map \([\blank]\) introduced in Definition~\ref{defn:grav} is additive on a disjoint union of constructible subsets,  take the union over the cover \(\{U\}\) of \(Y\) to see that \(\textstyle{\bigcap_{i = 1,\dots r}}F_{D_i}^{-1}(m_i) = \pi_t^{-1}(B_t)\) where
\[ [B_t] = \left[D_J^\circ\times\C^{tn - \sum _{j\in J} m_j}\times\left(\C^*\right)^{\vert J\vert}\right] = [D_J^\circ] \cdot\lm^{tn - \sum _{j\in J} m_j} \cdot (\lm - 1)^{\vert J\vert}.\]
Since \(\mu\left(\pi_t^{-1}(B_t)\right) = [B_t] \cdot \lm^{-(n + nt)} \),  we have
\[\mu\left(\textstyle{\bigcap_{i = 1,\dots r}}F_{D_i}^{-1}(m_i)\right) = [D_J^\circ ] \cdot \lm^{- \sum _{j\in J} m_j}\cdot (\lm - 1)^{\vert J\vert}\cdot \lm^{-n}.\]
Now use the partition (\ref{eqn:F_D_partition}) of \(F_D^{-1}(s)\) to compute the motivic integral:
\begin{eqnarray*}
\lefteqn{\sum_{s\in \mathbb{Z}_{\geq 0}} \mu\left(F_D^{-1}(s)\right) \cdot \lm^{-s} } \\
& = &  \sum_{s\in \mathbb{Z}_{\geq 0}} \sum_{J\subset \{1,\dots ,r\} } \sum_{(m_1,\dots ,m_r)\in M_{J,s}} \mu \left(\textstyle{\bigcap_{i = 1,\dots r}}F_{D_i}^{-1}(m_i)\right)\cdot \lm^{-\sum_{j\in J}a_j m_j} \\
 & = &  \sum_{s\in \mathbb{Z}_{\geq 0}}\sum_{J\subset \{1,\dots ,r\} } \sum_{(m_1,\dots ,m_r)\in M_{J,s}}[D_J^\circ ] \cdot (\lm - 1)^{\vert J\vert}\cdot \lm^{-n}\cdot \prod_{j\in J}\lm^{-(a_j + 1) m_j} \\
 & = & \sum_{J\subset \{1,\dots ,r\} }[D_J^\circ ] \cdot \prod_{j\in J} \left(\textstyle{ (\lm - 1)\cdot \sum_{m_j > 0} \lm^{-(a_j + 1) m_j} }\right) \cdot \lm^{-n} \\
 & = & \sum_{J\subset \{1,\dots ,r\} }[D_J^\circ ] \cdot \prod_{j\in J} \left(  (\lm - 1)\cdot \left(\frac{1}{1 - \lm^{-(a_j + 1)}} - 1\right)\right) \cdot \lm^{-n} \\
 & = & \sum_{J\subset \{1,\dots ,r\} }[D_J^\circ ] \cdot \left( \prod_{j\in J} \frac{ \lm - 1 }{\lm^{a_j + 1} - 1}\right) \cdot \lm^{-n}.
\end{eqnarray*}
\qed

 \begin{warning} 
 \label{warn:omission}
 \emph{There is a small error in the proof of the corresponding result in Batyrev~\cite[\S 6.28]{Batyrev:shn} which leads to the omission of the \(\lm^{-n}\) term.}
 \end{warning}

 \begin{corollary}
 \label{cor:subring} 
 The motivic integral of the pair \( (Y,D)\) is an element of the subring 
 \[
 \phi\bigl{(}K_0(\Var)[\lm^{-1}]\bigr{)}\left[\left\{ \frac{1}{\lm^i - 1}\right\}_{i\in\mathbb{N}}\right]
 \]
of the ring \(R\) introduced in Definition~\ref{defn:R}.
 \end{corollary}

 \subsection{The transformation rule for the integral}
The discrepancy divisor \(W := K_{Y'} - \alpha^*K_Y\) of a proper birational morphism \(\alpha \colon Y' \rightarrow Y\) between smooth varieties is the divisor of the Jacobian determinant of \(\alpha\).  The next result may therefore be viewed as the `change of variables formula' for the motivic integral.

 \begin{theorem}
 \label{thm:change_var} 
 Let \(\alpha \colon Y' \longrightarrow Y\) be a proper birational morphism of between smooth varieties and let \(W := K_{Y'} - \alpha^*K_Y\) be the discrepancy divisor.  Then 
 \[
 \int_{J_\infty(Y)} F_{D}\dmu = \int_{J_\infty(Y')}F_{\alpha ^*D + W} \dmu.
 \]
 \end{theorem}
 \proof Composition defines maps \(\alpha_t \colon J_t(Y') \rightarrow J_t(Y)\) for each \(t\in \mathbb{Z}_{\geq 0} \cup \oo\).  An arc in \(Y\) which is not contained in the locus of indeterminacy of \(\alpha^{-1}\) has a birational transform as an arc in \(Y'\).  In light of (\ref{eqn:zero}) and Proposition~\ref{prop:FDmble}, \(\alpha_\oo\) is bijective off a subset of measure zero.  

 \indent The sets \(F_W^{-1}(k)\),  for \(k\in \mathbb{Z}_{\geq 0}\),  partition \(J_\oo(Y') \setminus F_W^{-1}(\oo)\).  Thus,  for any \(s\in \mathbb{Z}_{\geq 0}\) we have,  modulo the set \(F_W^{-1}(\oo)\) of measure zero,  a partition
\begin{equation}
\label{eqn:C_k,s}
F_D^{-1}(s) = \bigsqcup_{k\in\mathbb{Z}_{\geq 0}} \alpha_\oo (C_{k,s}) \quad \mbox{where}\quad C_{k,s} := F_W^{-1}(k) \cap F_{\alpha^*D}^{-1}(s).
\end{equation}
The set \(C_{k,s}\) is cylinder and,  since the image of a constructible set is constructible (\cite[p.\ 72]{Mumford:rb}),  the set \(\alpha_\oo(C_{k,s})\) is cylinder.  Lemma~\ref{lemma:jacdet} below states that \(\mu\big{(}C_{k,s}\big{)} = \mu\big{(}\alpha_\oo(C_{k,s})\big{)}\cdot \lm^{k}\).  We use this identity and the partition (\ref{eqn:C_k,s}) to calculate
\[\int_{J_\infty(Y)} F_{D}\dmu = \sum_{k,s\in \mathbb{Z}_{\geq 0}} \mu\big{(}\alpha_\oo (C_{k,s})\big{)} \cdot \lm^{-s} = \sum_{k,s\in \mathbb{Z}_{\geq 0}} \mu\big{(}C_{k,s}\big{)}\cdot \lm^{-(s+k)}.\]
Set \(s' := s+k\).  Clearly \(\bigsqcup_{0\leq k\leq s'} C_{k,s'-k} = F_{\alpha^*D + W}^{-1}(s')\).  Substituting this into the above leaves
\[\int_{J_\infty(Y)} F_{D} \dmu = \sum_{s'\in\mathbb{Z}_{\geq 0}} \mu\big{(}F_{\alpha^*D + W}^{-1}(s') \big{)}\cdot \lm^{-s'} =  \int_{J_\infty(Y')} F_{\alpha ^*D + W} \dmu,\]
as required.\qed

 \begin{lemma} 
 \label{lemma:jacdet}
 \(\mu\big{(}C_{k,s}\big{)} = \mu\big{(}\alpha_\oo(C_{k,s})\big{)}\cdot \lm^{k}\).
 \end{lemma}
\discussproof
Both \(C_{k,s}\) and \(\alpha_\oo(C_{k,s})\) are cylinder sets so there exists \(t\in\mathbb{Z}_{\geq 0}\) and constructible sets \(B_t'\) and \(B_t\) in \(J_\oo(Y')\) and \(J_\oo(Y)\) respectively such that the following diagram commutes:
\[\begin{array}{ccc}
C_{k,s}\subset J_\oo(Y') & \stackrel{\alpha_\oo}{\longrightarrow} & \alpha_\oo(C_{k,s})\subset J_\oo(Y) \\
\pi_t\:\Big{\downarrow}\;  &  &  \;\Big{\downarrow}\:\pi_t \\
B_t'\subset J_t(Y') & \stackrel{\alpha_t}{\longrightarrow} & B_t \subset J_t(Y).
\end{array}\]  
We claim that the restriction of \(\alpha_t\) to \(B_t'\) is a \(\C^k\)-bundle over \(B_t\).  It follows that \([B_t'] = [\C^k]\cdot [B_t]\) and we have
\[\mu\big{(}C_{k,s}\big{)} = [B_t']\cdot \lm^{-(n+nt)} = [B_t]\cdot \lm^{k}\cdot\lm^{-(n+nt)} = \mu\big{(}\alpha_\oo(C_{k,s})\big{)}\cdot \lm^{k}\]
as required.  The proof of the claim is a local calculation which is carried out in \cite[Lemma 3.4(b)]{Denef:goa}.  The key observation is that the order of vanishing of the Jacobian determinant of \(\alpha\) at \(\gamma_y\in C_{k,s}\) is \(F_W(\gamma_y) = k\).
\qed

\begin{definition} \label{defn:misv} \emph{Let \(X\) denote a complex algebraic variety with at worst Gorenstein canonical singularities.  The \emph{motivic integral} of \(X\) is defined to be the motivic integral of the pair \( (Y,D)\),  where \(\varphi \colon Y \rightarrow X\) is any resolution of singularities for which the discrepancy divisor \(D = K_{Y} - \varphi^{*}K_{X}\) has only simple normal crossings.}
\end{definition}

Note first that the discrepancy divisor \(D\) is effective because \(X\) has at worst Gorenstein canonical singularities.  The crucial point however is that the motivic integral of \( (Y,D)\) is independent of the choice of resolution:

\begin{proposition}\label{prop:equalmotint} Let \(\varphi_1 \colon Y_1 \longrightarrow X\) and \(\varphi_2 \colon Y_2 \longrightarrow X\) be resolutions of \(X\) with discrepancy divisors \(D_1\) and \(D_2\) respectively.  Then the motivic integrals of the pairs \( (Y_1,D_1)\) and \( (Y_2,D_2)\) are equal.\end{proposition}

\proof Form a `Hironaka hut'
\[\begin{array}{rcl}
Y_{0} & \stackrel{\psi_2}{\longrightarrow} & Y_{2} \\
\psi_1\:\Big{\downarrow}\;  & \searrow &  \;\Big{\downarrow}\:\varphi_2\\
Y_1 & \stackrel{\varphi_1}{\longrightarrow} & X
\end{array}\]
and let \(D_0\) denote the discrepancy divisor of \(\varphi_0 \colon Y_0 \longrightarrow X\).  The discrepancy divisor of \(\psi_{i}\) is \(D_{0} - \psi_{i}^{*}D_{i}\).  Indeed
\begin{eqnarray*}
K_{Y_0} & = & \varphi_0^*(K_X) + D_0 = \psi_i^* \circ \varphi_i^*(K_X) + D_0  = \psi_i^* (K_{Y_i} - D_i) + D_0  \\
        & = & \psi_{i}^{*} (K_{Y_{i}}) + (D_{0} - \psi_{i}^{*}D_{i}).
\end{eqnarray*}
The maps \(\psi_i \colon Y_0 \longrightarrow Y_i\) are proper birational morphisms between smooth projective varieties so Theorem~\ref{thm:change_var} applies:
 \[
 \int_{J_\infty(Y_i)} F_{D_i} \dmu = \int_{J_\infty(Y_0)}F_{\psi_i^*D_i + (D_0 - \psi_i^*D_i)} \dmu = \int_{J_\infty(Y_0)} F_{D_0} \dmu.
 \]
This proves the result.\qed

 \section{Hodge numbers via motivic integration}
 \label{sec:hnvmi}
 This section describes how the motivic integral of \(X\) gives rise to the so-called ``stringy \(E\)-function'' which encodes the Hodge--Deligne numbers of a resolution \(Y\rightarrow X\). 

 \subsection{Encoding Hodge--Deligne numbers}
 \label{sec:ehdn}
 Deligne~\cite{Deligne:tdhII,Deligne:tdhIII} showed that the cohomology groups \(H^{k}(X,\Q)\) of a complex algebraic variety \(X\) carry a natural mixed Hodge structure.  This consists of an increasing weight filtration 
 \[
 0 = W_{-1} \subseteq W_{0} \subseteq \cdots \subseteq W_{2k} = H^{k}(X,\Q)
 \]
on the rational cohomology of \(X\) and a decreasing Hodge filtration
 \[
 H^{k}(X,\C) = F^{0}\supseteq F^{1} \supseteq \cdots\supseteq F^{k} \supseteq F^{k+1} = 0
 \]
on the complex cohomology of \(X\) such that the filtration induced by \(F^{\bullet}\) on the graded quotient \(\Gr_{l}^{W}H^{k}(X) := W_{l}/W_{l-1}\) is a pure Hodge structure of weight \(l\).  Thus
 \[
\Gr_{l}H^{k}(X)\otimes\:\C = F^{p}\Gr_{l}^{W}H^{k}(X)\oplus \overline{F^{l-p+1}\Gr_{l}^{W}H^{k}(X)}
 \]
 where \(F^{p}\Gr_{l}^{W}H^{k}(X)\) denotes the complexified image of \(F^{p}\cap W_{l}\) in the quotient \(W_{l}/W_{l-1}\otimes\:\C\).  The integers 
 \[
 h^{p,q}\big{(}H^{k}(X,\C)\big{)}:= \dimn_{\C} \left(F^{p}\Gr_{p+q}^{W}H^{k}(X)\cap \overline{F^{q}\Gr_{p+q}^{W}H^{k}(X)}\right)
 \]
are called the \emph{Hodge--Deligne numbers} of \(X\).  For a smooth projective variety \(X\) over \(\C\),  \(\Gr_{l}^{W}H^{k}(X,\Q) = 0\) unless \(l = k\) in which case the Hodge--Deligne numbers are the classical Hodge numbers \(h^{p,q}(X)\).

 \indent Danilov and Khovanskii~\cite{Danilov:hdn} observed that cohomology with compact support \(H_{c}^{k}(X,\Q)\) also admits a mixed Hodge structure and they encode the corresponding Hodge--Deligne numbers in a single polynomial:

 \begin{definition}
 \label{defn:E-poly}
 \emph{The \emph{\(E\)-polynomial} \(E(X)\in \Z[u,v]\) of a complex algebraic variety \(X\) of dimension \(n\) is defined to be 
 \[
 E(X) := \displaystyle{\sum_{0 \leq p,q \leq n} \sum_{0 \leq k\leq 2n}} (-1)^k h^{p,q}\left(H^k_c(X,\mathbb{C})\right) u^{p}\;v^{q}.
 \]
Evaluating \(E(X)\) at \(u = v = 1\) produces the standard topological Euler number \(e_{c}(X) = e(X)\).}
 \end{definition}
 
 \begin{theorem}[\cite{Danilov:hdn}]
 \label{thm:Epoly}
Let \(X, Y\) be complex algebraic varieties.  Then
 \begin{enumerate}
 \item[\one] if \(X = \bigsqcup X_{i}\) is stratified by a disjoint union of locally closed subvarieties then the \(E\)-polynomial is additive, i.e., \ \(E(X) = \sum E(X_{i})\).
 \item[\two] the \(E\)-polynomial is multiplicative,  i.e.,\ \(E(X\times Y) = E(X)\cdot E(Y)\).
 \item[\three] if \(f\colon Y \rightarrow X\) is a locally trivial fibration w.r.t.\ the Zariski topology and \(F\) is the fibre over a closed point then \(E(Y) = E(F)\cdot E(X)\).
 \end{enumerate}
 \end{theorem}

 See Danilov and Khovanskii~\cite{Danilov:hdn} for a proof.  

 \subsection{Kontsevich's theorem}
 \label{sec:kt}
 Theorem~\ref{thm:Epoly} asserts that the map \(E \colon \Var \longrightarrow \mathbb{Z}[u,v]\) associating to each complex variety \(X\) its \(E\)-polynomial is additive on a disjoint union of locally closed  subvarieties,  and satisfies \(E(X\times Y) = E(X)\cdot E(Y)\).  It follows from the universality of the map \([ \blank ]\) introduced in Definition~\ref{defn:grav} that \(E\) factors through the Grothendieck ring of algebraic varieties,  inducing a function \(E\colon K_0(\Var) \rightarrow  \mathbb{Z}[u,v]\).  By defining \(E(\lm^{-1}):= (uv)^{-1}\),  this extends to\footnote{One can use this function to define a finitely additive \(\mathbb{Z}[u,v,(uv)^{-1}]\)--valued measure \(\widetilde{\mu}_E := E\circ\widetilde{\mu}\) on cylinder sets given by \(\pi_k^{-1}(B_k) \rightarrow E(B_k)\cdot (uv)^{-n(k+1)}\).  Then construct the stringy \(E\)-function directly;  this is the approach adopted by Batyrev~\cite[\S 6]{Batyrev:shn}.}.
 \[
 E\colon K_0(\Var)[\lm^{-1}] \rightarrow  \mathbb{Z}[u,v,(uv)^{-1}].
 \]

 \begin{proposition} \label{prop:subring}
The map \(E\) can be extended uniquely to the subring \[\phi\bigl{(}K_0(\Var)[\lm^{-1}]\bigr{)}\left[\left\{ \frac{1}{\lm^i - 1}\right\}_{i\in\mathbb{N}}\right]\]
of the ring \(R\) introduced in Definition~\ref{defn:R}.
 \end{proposition}
 \proof  The kernel of the completion map \(\phi\colon K_0(\Var)[\lm^{-1}] \rightarrow R\) is 
 \begin{equation}
 \label{eqn:kernel}
 \bigcap_{m\in\mathbb{Z}} F^mK_0(\Var)[\lm^{-1}].
 \end{equation}  
For \([V]\cdot \lm^{-i}\in F^mK_0(\Var)[\lm^{-1}]\),  the degree of the \(E\)-polynomial \(E\left([V]\cdot \lm^{-i}\right)\) is \(2\:\mbox{dim}V - 2i \leq -2m\).  The \(E\)-polynomial of an element \(Z\) in the intersection (\ref{eqn:kernel}) must therefore be \( -\oo\);  that is,  \(E(Z) = 0\).  Thus \(E\) annihilates \(\mbox{ker }\phi\) and hence factors through \(\phi\left(K_0(\Var)[\lm^{-1}]\right)\).  The result follows when we define \(E(1/(\lm^i - 1)) := 1/((uv)^i - 1)\) for \(i\in \mathbb{N}\).\qed

\medskip

\para By Corollary~\ref{cor:subring} the motivic integral of the pair \( (Y,D)\) lies in the subring of Proposition~\ref{prop:subring}.  We now consider the image of the integral under \(E\).

 \begin{warning}
 \label{warning:error}
 \emph{As Warning~\ref{warn:omission} states,  the derivation of the motivic integral in \cite{Batyrev:shn} contains a small error which leads to the omission of an \(\lm^{-n}\) term.  However,  in practise it is extremely convenient to omit this term (!).  As a result,  we define the stringy \(E\)-function to be the image under \(E\) of the motivic integral times \(\lm^{n}\).  In short,  our stringy \(E\)-function agrees with that in \cite{Batyrev:shn},  even though our calculation of the motivic integral differs.}
 \end{warning}

 \begin{definition}
 \label{defn:mihd} 
 \emph{Let \(X\) be a complex algebraic variety of dimension \(n\) with at worst Gorenstein canonical singularities.  Let \(\varphi \colon Y \rightarrow X\) be a resolution of singularities for which the discrepancy divisor \(D = \sum_{i=1}^{r} a_iD_i\) has only simple normal crossings.  The \emph{stringy \(E\)-function} of \(X\) is 
 \begin{eqnarray}
E_{\mathrm{st}}(X) & := & E\left(\int_{J_\oo(Y)} F_{D} \dmu\cdot \:\lm^{n}\right) \nonumber \\
                   & = & \sum_{J\subseteq\{1,\dots ,r\} }E(D_J^\circ)\cdot\left(\prod_{j\in J}\frac{uv - 1}{(uv)^{a_j + 1} - 1}\right), \label{eqn:mihd}
 \end{eqnarray}
where we sum over all subsets \(J\subseteq\{1,\dots ,r\}\) including \(J = \emptyset\).}
 \end{definition}

\begin{theorem}[\cite{Kontsevich:lao}] \label{thm:crepant} Let \(X\) be a complex projective variety with at worst Gorenstein canonical singularities.  If \(X\) admits a crepant resolution \(\varphi \colon Y \rightarrow X\) then the Hodge numbers of \(Y\) are independent of the choice of crepant resolution.
\end{theorem}
\proof The discrepancy divisor \(D = \sum_{i=1}^r a_i D_i\) of the crepant resolution \(\varphi \colon Y \rightarrow X\) is by definition zero,  so the motivic integral of \(X\) is the motivic integral of the pair \( (Y,0)\).  Since each \(a_i = 0\) it's clear that
\[E_{\mathrm{st}}(X) = \textstyle{\sum_{J\subseteq\{1,\dots ,r\} }}E(D_J^\circ) = E(Y).\]
The stringy \(E\)-function is independent of the choice of the resolution \(\varphi\).  In particular,  \(E(Y) = E_{\mathrm{st}}(X) = E(Y_2)\) for \(\varphi_2 :\: Y_2 \rightarrow X\) another crepant resolution.  It remains to note that \(E(Y)\) determines the Hodge--Deligne numbers of \(Y\),  and hence the Hodge numbers since \(Y\) is smooth and projective.\qed

 \section{Calculating the motivic integral}
 \label{sec:ctmi}
 To perform nontrivial calculations of the stringy \(E\)-function we must choose varieties which admit no crepant resolution.  A nice family of examples is provided by Gorenstein terminal cyclic quotient singularities.

 \subsection{Toric construction of cyclic quotient singularities}
 \label{sec:aqs}
 Consider the action of the cyclic group \(G = \mathbb{Z}/r\subset \GL(n,\C)\) generated by the diagonal matrix\footnote{It is convenient to assume that \(\gcd (r,\alpha_{1},\dots ,\widehat{\alpha_{j}}, \dots ,\alpha_{n}) = 1\) for all \(j = 1,\dots ,n\) to ensure that the group action is `small'.  The notation \(\widehat{\alpha_{j}}\) means that \(\alpha_{j}\) is omitted.} 
 \[
 g = \diag\left(e^{2\pi i\alpha_{1}/r},\dots ,e^{2\pi i\alpha_{n}/r}\right) \quad \mbox{with } \quad 0\leq \alpha_{j} < r,
 \]
 where \(i = \sqrt{-1}\).  The quotient \(\C^{n}/G\) is the cyclic quotient singularity of type \(\frac{1}{r}(\alpha_{1},\dots ,\alpha_{n})\).  This fractional notation derives from the construction of \(\C^{n}/G\) as an affine toric variety as we now describe (see Reid~\cite[\S 4]{Reid:ypg} for more details).  

 Write \(\overline{M} \cong \Z^n\) for the lattice of Laurent monomials in \(x_{1}, \dots ,x_{n}\),  and \(\overline{N}\) for the dual lattice with basis \(e_{1},\dots ,e_{n}\).  Let \(\sigma = \R_{\geq 0}e_{1} + \dots + \R_{\geq 0}e_{n}\) denote the positive orthant in \(\overline{N}\otimes \R\) with dual cone \(\sigma^{\vee}\subset \overline{M}\otimes \R\).  The overlattice
 \begin{equation}
 \label{eqn:L}
N := \overline{N} + \mathbb{Z}\cdot \textstyle{\frac{1}{r}}(\alpha_{1},\dots ,\alpha_{n})
 \end{equation}
 is dual to \(M := \Hom(N,\Z)\).  A Laurent monomial in \(x_{1}, \dots ,x_{n}\) lies in the sublattice \(M\subset \overline{M}\) if and only if it is invariant under the action of the group \(G\).  Restricting to Laurent monomials with only nonnegative powers leads to the equality \(\C[x_{1},\dots ,x_{n}]^{G} = \C[\sigma^{\vee}\cap M]\), and hence
 \[
 \C^{n}/G = \Spec \:\C[x_{1},\dots ,x_{n}]^{G} = \Spec \:\C[\sigma^{\vee}\cap M] =: X_{N,\sigma}.
 \]
 In order to consider only Gorenstein terminal cyclic quotient singularities we impose certain restrictions on the type \(\frac{1}{r}(\alpha_{1},\dots ,\alpha_{n})\).  Watanabe~\cite{Watanabe:cis} showed that for a small subgroup \(G \subset \GL(n,\C)\),   \(\C^{n}/G\) is Gorenstein if and only if \(G \subset \SL(n,\C)\).  Thus, 
 \[
 \textstyle{X_{N,\sigma}\mbox{ is Gorenstein } \iff \sum_{j=1}^{n} \alpha_{j} \equiv 0 \mod{r}.}
 \]
 To determine when \(X_{N,\sigma}\) is terminal we recall the discrepancy calculation for cyclic quotients following Reid~\cite{Reid:c3f,Reid:ypg}.  Write \(\Box\) for the unit box in \(N\otimes\R\cong\R^{n}\), i.e., the unit cell of the sublattice \(\overline{N}\cong \Z^{n}\). Each element \(g \in G \cong N/\overline{N}\) has a unique representative 
 \[
 \textstyle{v_{g} = \frac{1}{r(g)}\big{(}\alpha_{1}(g),\dots ,\alpha_{n}(g)\big{)}\in N\cap \Box;}
 \]
 here \(v_{g}\) denotes both the vector in \(N\otimes \R\) and the lattice point in \(N\).  For each primitive vector \(v_{g}\in N\cap \Box\),  the barycentric subdivision of \(\sigma\) at \(v_{g}\) determines a toric blow-up \(\varphi\colon B \rightarrow A\) of the cyclic quotient \(A = X_{N,\sigma}\).  The exceptional divisor is \(D = X_{N(\tau),\mbox{\scriptsize{Star}}(\tau)}\)\footnote{The cones in the fan of \(B\) containing \(\tau\) as a face define a fan \(\mbox{Star}(\tau)\) in \(N(\tau)\otimes \R\) where \(N(\tau) := N /(\tau \cap N)\).  The toric variety \(X_{N(\tau),\mbox{\scriptsize{Star}}(\tau)}\) is the closure of the orbit \(O_{N,\tau}\).  See Fulton~\cite[p.~52]{Fulton:itv} for a nice picture.},  where \(\tau\) is the ray with primitive generator \(v_{g}\).  Adjunction for the toric blow-up \(\varphi\) is 
 \begin{equation}
 \label{eqn:discrepancy}
 K_{B} = \varphi^{*}K_{A} + \left(\frac{1}{r(g)}\sum_{j=1}^{n} \alpha_{j}(g) - 1\right)D,
 \end{equation}
 as shown by Reid~\cite[\S 4.8]{Reid:ypg}.  Thus, 
 \[
 X_{N,\sigma} \mbox{ is terminal }\iff \sum_{j=1}^{n} \frac{1}{r(g)}\alpha(g)_{j} > 1 \mbox{ for each }g\in G.
 \]
 That is,  every point \(v_{g}\in N\cap \Box\) must lie above the hyperplane \(\sum x_{i} = 1\).


  \para Given any fan \(\Sigma\) in the vector space \(N\otimes \R\),  the \(E\)-polynomial of the corresponding toric variety \(X_{N,\Sigma}\) is computed using the following simple formula:
 
 \begin{proposition}
 \label{prop:Epolytoric}
For a toric variety \(X_{N,\Sigma}\) of dimension \(n\) we have
 \begin{equation}
 \label{eqn:Epolytoric}
E\big{(}X_{N,\Sigma}\big{)} = \sum_{k=0}^{n} d_{k} \cdot \big{(}uv - 1\big{)}^{n-k},
 \end{equation}
where \(d_{k}\) is the number of cones of dimension \(k\) in \(\Sigma\).
 \end{proposition}
 \proof
 The Hodge numbers of \(\mathbb{P}^{1}\) are well known and,  by Theorem~\ref{thm:Epoly},  we compute \(E(\C^{*}) = E(\mathbb{P}^{1}) - E(\{0\}) - E(\{\oo\}) = (uv + 1) - 2 = uv - 1\).  The \(E\)-polynomial is multiplicative so \(E\big{(} (\C^{*})^{n-k}\big{)} =  E\big{(}\C^{*}\big{)}^{n-k} =\big{(}uv - 1\big{)}^{n-k}\). The action of the torus \(\mathbb{T}^{n}\simeq \left(\C^{*}\right)^{n}\) on \(X_{N,\Sigma}\) induces a stratification of \(X_{N,\Sigma}\) into orbits of the torus action \(O_{L,\tau} \cong \left(\C^{*}\right)^{n-\dimn \tau}\),  one for each cone \(\tau\in \Sigma\). The result follows from Theorem~\ref{thm:Epoly}\one.\qed

 \subsection{The examples}
 \label{sec:examples}
 For a finite subgroup \(G\subset \SL(n,\C)\) with \(n=2\) or \(3\),  the Gorenstein quotient \(\C^{n}/G\) admits a crepant resolution;  this is a classical result for surfaces and was established for 3-folds through a case by case analysis of the finite subgroups of \(\SL(3,\C)\) by Ito,  Markushevich and Roan (see Roan~\cite{Roan:mrg} and references therein).  We therefore begin by considering \hbox{4-dimensional} quotient singularities of type \(\frac{1}{r}(1,r-1,a,r-a)\) with \(\gcd(r,a) = 1\). Morrison and Stevens~\cite[Theorem 2.4(ii)]{Morrison:tqs} prove that these are the only Gorenstein terminal 4-fold cyclic quotient singularities.

 \begin{remark}
 \emph{In each example below we calculate both \(E(Y)\) and \(E_{\mathrm{st}}(\C^{n}/G)\) after resolving the singularity \(\varphi\colon Y\rightarrow \C^{n}/G\).  Note that \(E(Y)\) is not equal to the \(E\)-polynomial of the exceptional fibre \(D = \varphi^{-1}(\pi(0))\), for \(\pi\colon \C^{n}\rightarrow \C^{n}/G\) the quotient map.  Indeed,
 \[
 E(Y) = E(Y\setminus D) + E(D) = (uv)^{n} - 1 + E(D).
 \]
 The point is that the \(E\)-polynomial encodes the Hodge--Deligne numbers of compactly supported cohomology,  yet
\[
H^{*}_{c}(D,\C) = H^{*}(D,\C) \cong H^{*}(Y,\C) \neq H^{*}_{c}(Y,\C);
 \]
 the first equality holds because \(D\) is compact,  and the isomorphism is induced by a deformation retraction of \(Y\) onto \(D\subset Y\).}
 \end{remark}

 \begin{example}
 \label{ex:1/2(1111)}
 \emph{Write \(X = X_{N,\sigma}\) for the cyclic quotient singularity of type \(\frac{1}{2}(1,1,1,1)\).  Add the ray \(\tau\) generated by the vector \(v = \frac{1}{2}(1,1,1,1)\) to the cone \(\sigma\),  then take the simplicial subdivision of \(\sigma\).  This determines a toric resolution \(\varphi \colon Y \rightarrow X\) with a single exceptional divisor \(D = X_{N(\tau),\mbox{\scriptsize{Star}}(\tau)} \cong \mathbb{P}^3\). The discrepancy of \(D\) is \(1\) by (\ref{eqn:discrepancy}).  Using Proposition~\ref{prop:Epolytoric} we calculate
 \begin{eqnarray*}
E(Y) & = & E(Y \setminus D) + E(\mathbb{P}^3) \\
     & = & \big{(}(uv)^{4} - 1\big{)} + \big{(}(uv)^{3} + (uv)^{2} + uv + 1 \big{)}\\ 
     & = & (uv)^{4} + (uv)^{3} + (uv)^{2} + uv.
 \end{eqnarray*}
Compare this with the stringy \(E\)-function:
 \[
 E_{\mathrm{st}}(X) = E(Y \setminus D) + E(\mathbb{P}^3)\cdot \frac{uv - 1}{(uv)^2 - 1} = (uv)^{4} + (uv)^{2}.
 \]}
 \end{example}

 \smallskip

 \begin{example}
 \label{ex:1/3(1212)}
 \emph{Write \(X = X_{N,\sigma}\) for the cyclic quotient singularity of type \(\frac{1}{3}(1,2,1,2)\).  Add rays \(\tau_{1}\) and \(\tau_{2}\) generated by the vectors \(v_{1} = \frac{1}{3}(1,2,1,2)\) and \(v_{2} = \frac{1}{3}(2,1,2,1)\) respectively to the cone \(\sigma\),  then take the simplicial subdivision of \(\sigma\).  The resulting fan \(\Sigma\) is determined by its cross-section \(\Delta_{2} := \sigma\cap \left(\sum x_{i} = 2\right)\) illustrated in Figure~\ref{fig:3.1212}.}

 \begin{figure}[!ht]
 \begin{center}\input{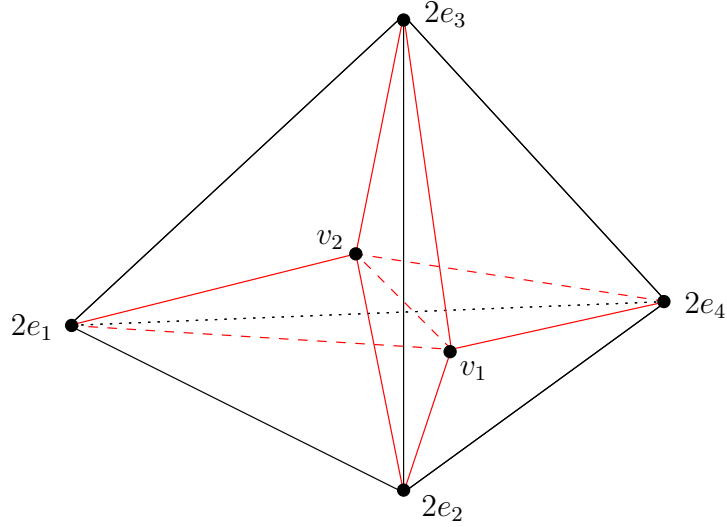}\end{center}
 \caption{The simplex $\Delta_{2}$ for $\frac{1}{3}(1,2,1,2)$}
 \label{fig:3.1212}
 \end{figure}

 \indent \emph{There are eight 3-dimensional simplices in \(\Delta_{2}\) (four contain a face of the tetrahedron and four contain the edge joining \(v_{1}\) to \(v_{2}\)).   Each of these simplices determines a 4-dimensional cone in \(\Sigma\) which is generated by a basis of the lattice \(N\),  so \(Y = X_{N,\Sigma}\to X_{N,\sigma}\) is a resolution.  The union of all eight 3-dimensional simplices in \(\Delta_{2}\) contain eighteen faces,  fifteen edges and six vertices.  Write \(d_{k}\) for the number of cones of dimension \(k\) in \(\Sigma\),  so
 \[
 d_{4} = 8;\quad d_{3} = 18;\quad d_{2} = 15;\quad d_{1} = 6;\quad d_{0} = 1 \;\mbox{(the origin in $N\otimes\R$)}.
 \]
Apply Proposition~\ref{prop:Epolytoric} to compute 
 \[
 E(Y) = (uv)^{4} + 2(uv)^{3} + 3(uv)^{2} + 2uv.
 \]
To compute \(E_{\mathrm{st}}(X)\) observe first that for \(j = 1\) or \(2\) the exceptional divisor \(D_{j} := X_{N(\tau_{j}),\mbox{\scriptsize{Star}}(\tau_{j})}\) has discrepancy \(1\) by (\ref{eqn:discrepancy}).  Write \(d_{k}(\tau_{j})\) for the number of cones of dimension \(k\) in \(\mbox{Star}(\tau_{j})\), so
 \[
 d_{3}(\tau_{j}) = 6;\quad d_{2}(\tau_{j}) = 9;\quad d_{1}(\tau_{j}) = 5;\quad d_{0}(\tau_{j}) = 1 \;\mbox{(the origin in \(N(\tau_{j})\))}.
 \]
Proposition~\ref{prop:Epolytoric} gives 
 \[
 E(D_{j}) = (uv)^{3} + 2(uv)^{2} + 2(uv) + 1\quad\mbox{for }j = 1,2.
 \]
 Similarly, the fan \(\mbox{Star}(\langle \tau_{1},\tau_{2}\rangle)\) contains four faces,  four edges and one vertex so Proposition~\ref{prop:Epolytoric} gives \(E(D_{1}\cap D_{2}) = (uv)^{2} + 2(uv) + 1\).  As a result
 \begin{itemize}
 \item \(E(D_{\emptyset}^{\circ}) = E(Y\setminus (D_{1}\cup D_{2})) = (uv)^{4} - 1\).
 \item \(E(D_{\{1\}}^{\circ}) = E(D_{\{2\}}^{\circ}) = E(D_{j}) - E(D_{1}\cap D_{2}) = (uv)^3 + (uv)^2\).
 \item \(E(D_{\{1,2\}}^{\circ}) = E(D_{1}\cap D_{2}) = (uv)^{2} + 2(uv) + 1\).
 \end{itemize}
 Now compute the stringy \(E\)-function using formula (\ref{eqn:mihd}):
 \begin{eqnarray*}
 E_{\mathrm{st}}(X) & = & (uv)^{4} - 1 + E(D_{\{1\}}^{\circ})\cdot \left(\frac{uv - 1}{(uv)^{2} - 1}\right) + E(D_{\{2\}}^{\circ})\cdot \left(\frac{uv - 1}{(uv)^{2} - 1}\right) \\
 & & \quad +\; E(D_{\{1,2\}}^{\circ})\cdot  \left(\frac{uv - 1}{(uv)^{2} - 1}\right)^{2} \\
 & = & (uv)^{4} + 2(uv)^{2}.
 \end{eqnarray*}
 }
 \end{example}

\begin{example}
\label{ex:1/4(1313)}
\emph{Write \(X = X_{N,\sigma}\) for the cyclic quotient singularity of type \(\frac{1}{4}(1,3,1,3)\).  Add rays \(\tau_{1}\), \(\tau_{2}\) and \(\tau_{3}\) generated by the vectors \(v_{1} = \frac{1}{4}(1,3,1,3)\) and \(v_{2} = \frac{1}{4}(2,2,2,2)\) and \(v_{3} = \frac{1}{4}(3,1,3,1)\) respectively to the cone \(\sigma\),  then take the simplicial subdivision of \(\sigma\).  The cross-section \(\Delta_{2}\) of the resulting fan \(\Sigma\) has three colinear points in the interior of the tetrahedron but is otherwise similar to that shown in Figure~\ref{fig:3.1212}.  There are twelve 3-dimensional simplices in \(\Delta_{2}\) containing 26 faces,  20 edges and 7 vertices.  Proposition~\ref{prop:Epolytoric} calculates
 \[
 E(Y) = (uv)^{4} + 3(uv)^{3} + 5(uv)^{2} + 3uv.
 \] 
 For \(j = 1,2,3\),  the divisors \(D_{j} := X_{N(\tau_{j}),\mbox{\scriptsize{Star}}(\tau_{j})}\) have discrepancy \(1\) by (\ref{eqn:discrepancy}).  Following the method of Example~\ref{ex:1/3(1212)} we calculate
 \[
 E(D_{1}) = E(D_{3}) = (uv)^{3} + 2(uv)^{2} + 2(uv) + 1
 \]
and \(E(D_{1}\cap D_{2}) = E(D_{2}\cap D_{3}) = (uv)^{2} + 2(uv) + 1\).  To compute the \(E\)-polynomial of \(D_{2}\) observe that
 \[
 d_{3}(\tau_{2}) = 8;\quad d_{2}(\tau_{2}) = 12;\quad d_{1}(\tau_{2}) = 6;\quad d_{0}(\tau_{2}) = 1 \;\mbox{(the origin in \(N(\tau_{2})\))},
 \]
 where \(d_{k}(\tau_{2})\) denotes the number of cones of dimension \(k\) in \(\mbox{Star}(\tau_{2})\).  It follows from Proposition~\ref{prop:Epolytoric} that
 \[
 E(D_{2}) = (uv)^{3} + 3(uv)^{2} + 3(uv) + 1.
 \]
 Finally,  since \(D_1 \cap D_3 = \emptyset\) we have \(E(D_1\cap D_3) = E(D_{1}\cap D_2\cap D_3) = 0\).  As a result
 \begin{itemize}
 \item \(E(D_{\emptyset}^{\circ}) = E(Y\setminus (D_{1}\cup D_{2}\cup D_{3})) = (uv)^{4} - 1\).
 \item \(E(D_{\{1\}}^{\circ}) = E(D_{\{3\}}^{\circ}) = (uv)^3 + (uv)^2\).
 \item \(E(D_{\{2\}}^{\circ}) = (uv)^3 + (uv)^2 - (uv) - 1\).
 \item \(E(D_{\{1,2\}}^{\circ}) = E(D_{\{2,3\}}^{\circ}) = (uv)^{2} + 2(uv) + 1\).
 \item \(E(D_{\{1,3\}}^{\circ}) = E(D_{\{1,2,3\}}^{\circ}) = 0\).
 \end{itemize}
 Apply formula (\ref{eqn:mihd}) to compute \(E_{\mathrm{st}}(X) = (uv)^{4} + 3(uv)^{2}\).}
 \end{example}

 \begin{remark} \emph{The above examples feature only exceptional divisors with discrepancy 1.  To obtain examples of Gorenstein terminal cyclic quotient singularities which admit resolutions containing divisors having discrepancy larger than one we must work in dimension higher than four.}
 \end{remark}

\begin{example}
\label{ex:1/r(111...1)}
\emph{Write \(X = X_{N,\sigma}\) for the cyclic quotient singularity of type \(\frac{1}{r}(1,1,1,\dots,1)\) where \(n := \dimn X = kr\) for some \(k\in \Z\) (by assuming that \(r\) divides \(n\) we ensure that \(X\) is Gorenstein).  Add a single ray \(\tau\) generated by the vector \(v_{1} = \frac{1}{r}(1,1,1,\dots,1)\) to the cone \(\sigma\),  then take the simplicial subdivision of \(\sigma\).  This determines a toric resolution \(\varphi \colon Y \rightarrow X\) with a single exceptional divisor \(D = X_{N(\tau),\mbox{\scriptsize{Star}}(\tau)} \cong \mathbb{P}^{n-1}\).  The discrepancy of \(D\) is \(k-1\) by (\ref{eqn:discrepancy}).  Using Proposition~\ref{prop:Epolytoric} we calculate
 \begin{eqnarray*}
E(Y) & = & E(Y \setminus D) + E(\mathbb{P}^{n-1}) \\
     & = & \big{(}(uv)^{n} - 1\big{)} + \big{(}(uv)^{n-1} + (uv)^{n-2} + \dots + uv + 1 \big{)}\\ 
     & = & (uv)^{n} + (uv)^{n-1} + \dots + uv.
 \end{eqnarray*}
Compare this with the stringy \(E\)-function:
 \begin{eqnarray*}
 E_{\mathrm{st}}(X) & = & E(Y \setminus D) + E(\mathbb{P}^{n-1})\cdot \frac{uv - 1}{(uv)^{k} - 1}\\ 
		& = & (uv)^{n} + (uv)^{n-k} + \dots + (uv)^{2k} + (uv)^{k}.
 \end{eqnarray*}
 }
 \end{example}

 \begin{example}
 \label{ex:1/3(121212)}
 \emph{Let \(X = X_{N,\sigma}\) denote the cyclic quotient singularity of type \(\frac{1}{3}(1,2,1,2,1,2)\) (compare Example~\ref{ex:1/3(1212)}).  Add rays \(\tau_{1}\) and \(\tau_{2}\) generated by the vectors \(v_{1} = \frac{1}{3}(1,2,1,2,1,2)\) and \(v_{2} = \frac{1}{3}(2,1,2,1,2,1)\) respectively to the cone \(\sigma\),  then take the simplicial subdivision of \(\sigma\).  Both \(v_{1}\) and \(v_{2}\) lie in the simplex \(\Delta_{3}:= \sigma\cap \left(\sum x_{i} = 0\right)\) of the resulting fan \(\Sigma\) so the corresponding exceptional divisors \(D_{1}\) and \(D_{2}\) each have discrepancy \(2\) by (\ref{eqn:discrepancy}).  The cross-section \(\Delta_{3}\) is difficult to draw (it is 5-dimensional!) but,  using Figure~\ref{fig:3.1212} as a guide,  one can show that 
 \[
 d_{6} = 15;\quad d_{5} = 48;\quad d_{4} = 68;\quad d_{3} = 56;\quad d_{2} = 28;\quad d_{1} = 8;\quad d_{0} = 1,
 \]
 where \(d_{k}\) denotes the number of cones of dimension \(k\) in \(\Sigma\).  Hence
 \[
 E(Y) = (uv)^{6} + 2(uv)^{5} + 3(uv)^{4} + 4(uv)^{3} + 3(uv)^{2} + 2(uv).
 \]
 As with Example~\ref{ex:1/3(1212)},  for \(j=1,2\) write \(d_{k}(\tau_{j})\) for the number of cones of dimension \(k\) in \(\mbox{Star}(\tau_{j})\), so
 \[
 d_{5}(\tau_{j}) = 12;\;d_{4}(\tau_{j}) = 30;\;d_{3}(\tau_{j}) = 34;\; d_{2}(\tau_{j}) = 21;\; d_{1}(\tau_{j}) = 7;\; d_{0}(\tau_{j}) = 1.
 \]
Proposition~\ref{prop:Epolytoric} gives 
 \[
 E(D_{j}) = (uv)^{5} + 2(uv)^{4} + 3(uv)^{3} + 3(uv)^{2} + 2(uv) + 1\quad\mbox{for }j = 1,2.
 \]
 Similarly,  counting simplices in the fan \(\mbox{Star}(\langle \tau_{1},\tau_{2}\rangle)\) gives
 \[
 E(D_{1}\cap D_{2}) = (uv)^{4} + 2(uv)^{3} + 3(uv)^{2} + 2(uv) + 1.
 \]
 Now compute the stringy \(E\)-function using formula (\ref{eqn:mihd}):
 \begin{eqnarray*}
 E_{\mathrm{st}}(X) & = & (uv)^{6} - 1 + E(D_{\{1\}}^{\circ})\cdot \left(\frac{uv - 1}{(uv)^{3} - 1}\right) + E(D_{\{2\}}^{\circ})\cdot \left(\frac{uv - 1}{(uv)^{3} - 1}\right) \\
 & & \quad +\; E(D_{\{1,2\}}^{\circ})\cdot  \left(\frac{uv - 1}{(uv)^{3} - 1}\right)^{2} \\
 & = & (uv)^{6} + 2(uv)^{3}.
 \end{eqnarray*}
 }
 \end{example}

 \section{The McKay correspondence}
 The stringy \(E\)-function of a Gorenstein canonical quotient singularity \(\C^{n}/G\) can be calculated in terms of the representation theory of the finite subgroup \(G\subset \SL(n,\C)\) using a simple formula due to Batyrev~\cite{Batyrev:nai,Batyrev:ca} (see also Denef and Loeser~\cite{Denef:mc}).  To state the formula,  note that each \(g\in G\) is conjugate to a diagonal matrix 
 \begin{equation}
 \label{eqn:diagmatrix}
 g = \diag\left(e^{2\pi i\alpha_{1}(g)/r(g)},\dots ,e^{2\pi i\alpha_{n}(g)/r(g)}\right) \quad \mbox{with } \quad 0\leq \alpha_{j}(g) < r(g),
 \end{equation}
where \(r(g)\) is the order of \(g\) and \(i = \sqrt{-1}\).  To each conjugacy class \( [g]\) of the group \(G\) we associate an integer in the range \(0\leq \age[g] \leq n-1\) defined by
 \[
 \age [g] := \frac{1}{r(g)}\sum_{j=1}^{n} \alpha_{j}(g).
 \]   
 (In particular,  for the cyclic action introduced in \S\ref{sec:aqs},  the age grading on \(G\) corresponds to the slicing of the unit box \(\Box \subset \overline{N}_{\R}\cong\R^{n}\) into polytopes \(\Delta_{k} := \sigma\cap \left(\sum x_{i} = k\right)\) for \(k = 0,\dots ,n-1\).) 

Batyrev's formula is 
 \begin{equation}
 \label{eqn:stringformula}
E_{\mathrm{st}}(\C^{n}/G) = \sum_{[g]\in \Conj(G)} (uv)^{n-\age[g]},
 \end{equation}
where we sum over the conjugacy classes of \(G\).  

 For example,  the nontrivial element \(g\) of the group \(G = \Z/2\) acting on \(\C^{4}\) in Example~\ref{ex:1/2(1111)} has \(\age\) two because \(v_{g} = \frac{1}{2}(1,1,1,1)\in \Delta_{2}\).  Formula (\ref{eqn:stringformula}) gives \(E_{\mathrm{st}}(\C^{4}/G) = (uv)^{4} + (uv)^{2}\) as shown in \S\ref{sec:examples}.  (The reader should check that the same holds for the other examples of \S\ref{sec:examples}.)  In this section we recall formula (\ref{eqn:stringformula}) in the wider context of the cohomological McKay correspondence.

 \subsection{Reid's McKay correspondence conjecture}
 \label{sec:rmc}
 Motivated by string theory,  Dixon et al.~\cite{Dixon:sooI} introduced the \emph{orbifold Euler number} for a finite group \(G\) acting on a manifold \(M\).  This number can be written in the form\footnote{The formula given here is due to Hirzebruch and H\"{o}fer~\cite{Hirzebruch:eno}.}
 \[
 e(M,G) = \sum_{[g]\in \Conj(G)} e\big{(}M^{g}/C(g)\big{)}, 
 \]
 where the sum runs over the conjugacy classes of \(G\), \(e\) denotes the topological Euler number,  \(M^{g}\) is the fixed point set of \(g\) and \(C(g)\) is the centraliser of \(g\).  Dixon et al.~\cite{Dixon:sooII} formulated what became known as the ``physicists' Euler number conjecture'':

 \begin{conjecture}
 \label{conj:penc}
 If \(M/G\) is a Gorenstein Calabi--Yau variety which admits a crepant resolution \(Y \rightarrow M/G\) then \(e(Y) = e(M,G)\).
 \end{conjecture}

 Hirzebruch and H\"{o}fer~\cite{Hirzebruch:eno} observed that for a finite subgroup \(G \subset \U(n)\) acting on \(M = \C^{n}\),  the orbifold Euler number is equal to the number of conjugacy classes of \(G\) because every fixed point set \(M^{g}\) is contractible.  For \(G \subset \SU(2,\C)\),  the classical McKay correspondence states that the second Betti number \(b_{2}(Y)\) of the minimal resolution \(Y \rightarrow \C^{2}/G\) equals the number of nontrivial conjugacy classes of \(G\).  As a result,  the equality
 \begin{equation}
 \label{eqn:HHMcKay}
 e(Y) = b_{2}(Y) + 1 = \#\Big{\{}\mbox{conjugacy classes of }G\Big{\}} = e(\C^{2},G)
 \end{equation}
can be viewed as a version of the McKay correspondence. Inspired by this observation,  Reid~\cite{Reid:pen} proposed that a version of the statement (\ref{eqn:HHMcKay}) should hold in arbitrary dimension:

 \begin{conjecture}[McKay correspondence] 
 \label{conj:generalised}
For \(G\subset \SL(n,\C)\) a finite subgroup,  suppose that the quotient variety \(X := \C^{n}/G\) admits a crepant resolution \(\varphi\colon Y\rightarrow X\).  Then \(H^{*}(Y,\Q)\) has a basis consisting of algebraic cycles corresponding one-to-one with conjugacy classes of \(G\).  In particular
 \[
 e(Y) = \#\Big{\{}\mbox{conjugacy classes of }G\Big{\}} = e(\C^{n},G).
 \]
 \end{conjecture}

 \subsection{McKay correspondence via motivic integration}
 \label{sec:mcvmi}
 The key to Batyrev's proof of Conjecture~\ref{conj:generalised} is formula (\ref{eqn:stringformula}):

 \begin{theorem}[\cite{Batyrev:nai}]
 \label{thm:stringy}
 Let \(G\subset \SL(n,\C)\) be a finite subgroup.  Then
 \[
 E_{\mathrm{st}}(\C^{n}/G) = \sum_{[g]\in \Conj(G)} (uv)^{n-\age[g]},
 \]
 where the sum runs over conjugacy classes of \(G\).
 \end{theorem}

 \noindent\textsc{Proof of the Abelian case\footnote{The proof of the general case is a consequence of the equality of the stringy \(E\)-function and the orbifold \(E\)-function introduced by Batyrev~\cite{Batyrev:nai}.}.\ } Choose coordinates on \(\C^{n}\) so that every matrix \(g\in G\) takes the form given in (\ref{eqn:diagmatrix}).  The construction of \S\ref{sec:aqs} can be adapted to show that \(\C^{n}/G\) is the toric variety \(X_{N,\sigma}\) corresponding to the cone \(\sigma = \R_{\geq 0}e_{1} + \dots + \R_{\geq 0}e_{n}\) and the lattice
 \begin{equation}
 \label{eqn:N}
N := \Z e_{1} + \dots + \Z e_{n} + \sum_{g\in G} \; \Z\cdot v_{g},  \quad \mbox{for } v_{g} = \frac{1}{r(g)}\big{(}\alpha_{1}(g),\dots ,\alpha_{n}(g)\big{)}.
 \end{equation}
 Now \(X_{N,\sigma}\) is Gorenstein because \(G\subset \SL(n,\C)\),  so there exists a continuous linear function \(\psi_{K}\colon N\otimes \R\to \R_{\geq 0}\) satisfying \(\psi(e_{i}) = 1\) for \(i = 1,\dots ,n\).  A straightforward computation due to Batyrev~\cite[Theorem 4.3]{Batyrev:shn} gives
 \begin{equation}
 \label{eqn:toricE}
 E_{\mathrm{st}}(X_{N,\sigma}) = (uv - 1)^{n} \sum_{v\in N\cap \sigma} (uv)^{-\psi_{K}(v)}.
 \end{equation}
 For \(v_{g}\in \Box\),  \(\psi_{K}(v_{g}) = \age(g)\).  In fact,  for any lattice point \(v\in N\cap \sigma\) we have \(\psi_{K}(v) = k \iff v\in\Delta_{k} = \sigma\cap \left(\sum x_{i} = k\right)\).  Thus for each \(v\in N\cap \sigma\) there exists unique \(v_{g}\in \Box\) and \((x_{1},\dots ,x_{n})\in \Z^{n}_{\geq 0}\) such that \(v\) is the translate of \(v_{g}\) by \((x_{1},\dots ,x_{n})\), and \(\psi_{K}(v) = \age(g) + \sum_{i = 1}^{n} x_{i}\). 

 Covering the positive orthant \(\sigma\) by translations of \(\Box\) gives
 \begin{eqnarray*}
 \sum_{v\in N\cap \sigma}(uv)^{-\psi_{K}(v)} & = & \sum_{v_{g}\in \Box}\sum_{(x_{1},\dots ,x_{n}) \in \Z^{n}_{\geq 0}} (uv)^{-\age(g) - \sum x_{i}} \\
 & = & \sum_{v_{g}\in \Box} (uv)^{-\age (g)} \prod_{i = 1}^{n} \sum_{x_{i} \in \Z_{\geq 0}} (uv)^{-x_{i}} \\
 & = & \sum_{g\in G} (uv)^{-\age (g)}\prod_{i = 1}^{n} \frac{1}{1 - (uv)^{-1}}.
 \end{eqnarray*}
 Substituting this into (\ref{eqn:toricE}) gives
  \[ 
 E_{\mathrm{st}}(X_{N,\sigma}) = \sum_{g\in G} (uv)^{-\age (g)}\prod_{i = 1}^{n} \frac{uv - 1}{1 - (uv)^{-1}} = \sum_{g\in G} (uv)^{n - \age (g)}.
 \]
 This proves the theorem for a finite Abelian subgroup \(G \subset \SL(n,\C)\). 
 \qed

 \begin{corollary}[strong McKay correspondence]
 \label{thm:smc}
 Let \(G\subset \SL(n,\C)\) be a finite subgroup and suppose that the quotient \(X = \C^{n}/G\) admits a crepant resolution \(\varphi\colon Y\rightarrow X\).  The nonzero Betti numbers of \(Y\) are
 \[
 \dimn_{\C} H^{2k}(Y,\C) = \#\Big{\{}\mbox{age } k \mbox{ conjugacy classes of }G\Big{\}}.
 \]
for \(k = 0,\dots,n-1\).  In particular,  Theorem~\ref{introthm:gmc} holds.
 \end{corollary}
 \proof  The Hodge structure in \(H^{i}_{c}(Y,\Q)\) is pure for each \(i\) and Poincar\'{e} duality \(H^{2n - i}_{c}(Y,\C) \otimes H^{i}(Y,\C) \to H^{2n}_{c}(Y,\C)\) respects the Hodge structure,  so it is enough to show that the only nonzero Hodge--Deligne numbers of the compactly supported cohomology of \(Y\) are
 \[
 h^{n-k,n-k}\big{(}H_{c}^{2n-2k}(Y,\C)\big{)} = \#\Big{\{}\mbox{age } k \mbox{ conjugacy classes of }G\Big{\}}.
 \]
 Now \(h^{n-k,n-k}\big{(}H_{c}^{2n-2k}(Y,\C)\big{)}\) is the coefficient of \((uv)^{n-k}\) in the \hbox{\(E\)-polynomial} of \(Y\).  Moreover,  the resolution \(\varphi\colon Y \to X\) is crepant so \(E(Y) = E_{\mathrm{st}}(X)\) and the result follows from Theorem~\ref{thm:stringy}.
 \qed

 \appendix

 \section{Why `motivic' integration?}
 \label{sec:mni}
 In this appendix we investigate the motivic nature of the integral.  We also justify the notation \lm\ for the class of the complex line \(\C\) in the Grothendieck ring of algebraic varieties.

\para The category \(\mathcal{M}_{\C}\) of Chow motives over \C\ is defined as follows (see \cite{Scholl:cm}): an object is a triple \( (X,p,m)\) where \(X\) is a smooth,  complex projective variety of dimension \(d\),  \(p\) is an element of the Chow ring \(A^d(X\times X)\) which satisfies \(p^2 = p\) and \(m\in \mathbb{Z}\).  If \( (X,p,m)\) and \( (Y,q,n)\) are motives then 
\[\Hom_{\mathcal{M}_{\C}}\left( (X,p,m),(Y,q,n)\right) = q A^{d + n - m}(X,Y)p\]
where composition of morphisms is given by composition of correspondences.  \(\mathcal{M}_{\C}\) is additive,  \(\mathbb{Q}\)-linear and pseudo-abelian.  Tensor product of motives is defined as \((X,p,m)\otimes (Y,q,n) = (X\times Y, p\otimes q, m+n)\).  There is a functor 
\[h\colon \Var^{\circ} \rightarrow \mathcal{M}_{\C}\] 
which sends \(X\) to \( (X,\Delta_X,0)\),  the Chow motive of \(X\),  where the diagonal \(\Delta_X \subset X\times X\) is the identity in \(A^*(X\times X)\).  The motive of a point \(1 = h(\Spec\ \C)\) is the identity with respect to tensor product.  The Lefschetz motive \(\lm\) is defined implicitly via the relation \(h(\mathbb{P}^1_{\C}) = 1 \oplus \lm\).
 
\begin{definition} \emph{The Grothendieck group of \(\mathcal{M}_{\C}\) is the free abelian group generated by isomorphism classes of objects in \(\mathcal{M}_{\C}\) modulo the subgroup generated by elements of the form \( [(X,p,m)] - [(Y,q,n)] - [(Z,r,k)]\) whenever \( (X,p,m) \simeq (Y,q,n) \oplus (Z,r,k)\).  Tensor product of motives induces a ring structure and the resulting ring,  denoted \(K_0(\mathcal{M}_{\C})\),  is the Grothendieck ring of Chow Motives (over \C).}\end{definition}

Gillet and Soul\'{e} \cite{Gillet:dmK} exhibit a map 
\[M\colon \Var \longrightarrow K_0(\mathcal{M}_{\C})\]
which sends a smooth, projective variety \(X\) to the class \( [h(X)]\) of the motive of \(X\).  Furthermore the map is additive on disjoint unions of locally closed subsets and satisfies \(M(X\times Y) = M(X)\cdot M(Y)\).  

\para We now play the same game as we did in \S\ref{sec:hnvmi}.  Namely,  \(M\) factors through \(K_0(\Var)\) inducing 
\[M \colon K_0(\Var) \longrightarrow K_0(\mathcal{M}_{\C}).\]  
Observe that the image of \([\C]\) under \(M\) is the class of the Lefschetz motive \(\lm\);  this explains why we use the notation \(\lm\) to denote the class of \(\C\) in \(K_0(\Var)\) in \S~\ref{sec:cmi}.  Sending \(\lm^{-1} \in K_0(\Var)[\lm^{-1}]\) to \(\lm^{-1} \in K_0(\mathcal{M}_{\C})\) produces a map  
\[M\colon K_0(\Var)[\lm^{-1}] \longrightarrow K_0(\mathcal{M}_{\C}).\]
At present it is unknown whether or not \(M\) annihilates the kernel of the natural completion map \(\phi \colon K_0(\Var)[\lm^{-1}] \rightarrow R\).  Denef and Loeser conjecture that it does (see \cite[Remark 1.2.3]{Denef:miz}).  If this is true,  extend \(M\) to a ring homomorphism 
 \[
 M_{\mathrm{st}}\colon \phi\left(K_0(\Var)[\lm^{-1}]\right)\left[\left\{\frac{1}{\lm^i - 1}\right\}_{i\in\mathbb{N}}\right]\longrightarrow K_0(\mathcal{M}_{\C})\left[\left\{\frac{1}{\lm^i - 1}\right\}_{i\in\mathbb{N}}\right]
 \]
such that the image of \([D_J^\circ]\) under \(M_{\mathrm{st}}\) is equal to \(M(D_J^\circ)\).

\begin{definition} \emph{Let \(X\) denote a complex algebraic variety with at worst canonical,  Gorenstein singularities and let \(\varphi \colon Y \rightarrow X\) be any resolution of singularities for which the discrepancy divisor \(D = \sum a_iD_i\) has only simple normal crossings.  The \emph{stringy motive} of \(X\) is
 \begin{eqnarray*}
 M_{\mathrm{st}}(X) & := & M_{\mathrm{st}}\left(\int_{J_\infty(Y)} F_{D}\dmu\cdot \; \lm^{n}\right) \\
           & =  & \sum_{J\subseteq\{1,\dots ,r\} }M(D_J^\circ)\cdot\left(\prod_{j\in J}\frac{\lm - 1}{\lm^{a_j + 1} - 1}\right)
 \end{eqnarray*}
where we sum over all subsets \(J\subseteq\{1,\dots ,r\}\) including \(J = \emptyset\).  As with the definition of the stringy \(E\)-function (see Definition~\ref{defn:mihd}) we multiply by \(\lm^{n}\) for convenience.}\end{definition}

 \hrulefill

 \bibliography{research}

 \medskip

 \noindent \textsc{Department of Mathematics, \\ 
 University of Utah, \\
 155 South 1400 East, \\
 Salt Lake City,
 UT 84112, USA. \\
 E-mail:} \texttt{craw@math.utah.edu}

 \end{document}